\documentclass{ifacconf}

\usepackage{graphicx}      
\usepackage{natbib}        

\usepackage[T1]{fontenc} 

\usepackage[utf8]{inputenc} 

\usepackage{graphicx} 
\graphicspath{{Figures/}} 

\usepackage{enumitem} 
\usepackage{array,booktabs}
\usepackage{subfig}
\usepackage{makecell}
\usepackage{float}
\usepackage{subfig} 

\usepackage{amsmath,amssymb} 
\usepackage{mathtools}
\usepackage{amsfonts}
\usepackage{bbm}
\usepackage{multirow}
\usepackage{lastpage}
\usepackage{scrtime}
\usepackage{mathrsfs,dsfont} 
\usepackage{eurosym}
\usepackage{amscd}

\usepackage{varioref} 
\usepackage{tikz}
\usetikzlibrary {positioning}
\usepackage{pgfplots}
\pgfplotsset{ 
	compat=newest, 
	legend style =
	{font=\small \sffamily},
	label style = {font=\small\sffamily},
	every tick label/.append style={font=\small}}
\newcommand{\mc}{\mathcal}
\usetikzlibrary{arrows}

\newcommand{\overbar}[1]{\mkern 1.5mu\overline{\mkern-1.5mu#1\mkern-1.5mu}\mkern 1.5mu}
\newcommand{\underbarr}[1]{\mkern 1.5mu\underline{\mkern-1.5mu#1\mkern-1.5mu}\mkern 1.5mu}

\def\qed{\hfill \vrule height 7pt width 7pt depth 0pt\medskip}
\newcommand{\ds}{\displaystyle}

\newcommand{\ba}{\begin{array}}
\newcommand{\ea}{\end{array}}

\newcommand{\1}{\mathbbm{1}}
\renewcommand{\l}{\left}\renewcommand{\r}{\right}
\newcommand{\be}{\begin{equation}}
\newcommand{\ee}{\end{equation}}
\newcommand{\eps}{\varepsilon}

\renewcommand{\mc}{\mathcal}

\newcommand{\ov}{\overline}
\newcommand{\ul}{\underline}

\newcommand{\R}{\mathbb{R}}

\def\R{\mathbb{R}}

\newtheorem{lemma}{Lemma}
\newtheorem{proposition}{Proposition}

\newtheorem{theorem}{Theorem}

\usetikzlibrary[topaths]
\newcount\mycount

\usepackage[framemethod=tikz]{mdframed}

\newcounter{examplecounter}

\definecolor{mycolor}{rgb}{0.122, 0.435, 0.698}
\tikzset{every picture/.style={line width=0.75pt}} 

\begin{document}
\begin{frontmatter}
	\title{Stability and phase transitions of dynamical flow networks with finite capacities}
	
	\thanks[footnoteinfo]{Giacomo Como is also with the Department of Automatic Control, Lund University, Sweden. This work was partially supported by MIUR grant Dipartimenti di Eccellenza 2018--2022 [CUP: E11G18000350001], the Swedish Research Council, and by the Compagnia di San Paolo.}
	
	\author[First]{Leonardo Massai} 
	\author[First]{Giacomo Como} 
	\author[First]{Fabio Fagnani}
	
	\address[First]{Department of Mathematical Sciences ``G.L.~Lagrange'', Politecnico di Torino, Corso Duca degli Abruzzi 24, 10129 Torino, Italy\\ 
 (e-mail: \{leonardo.massai,giacomo.como,fabio.fagnani\}@polito.it).}

\begin{abstract}                
We study deterministic continuous-time lossy dynamical flow networks with constant exogenous demands, fixed routing, and finite flow and buffer capacities. In the considered model, when the total net flow in a cell ---consisting of the difference between the total flow directed towards it minus the outflow from it--- exceeds a certain capacity constraint, then the exceeding part of it leaks out of the system. The ensuing network flow dynamics is a linear saturated system with compact state space that we analyse using tools from monotone systems and contraction theory. Specifically, we prove that there exists a set of equilibria that is globally asymptotically stable. Such equilibrium set reduces to a single globally asymptotically stable equilibrium for generic exogenous demand vectors. 
Moreover, we show that the critical exogenous demand vectors giving rise to non-unique equilibria correspond to phase transitions in the asymptotic behavior of the dynamical flow network. 
\end{abstract}

\begin{keyword}
Dynamical flow networks, nonlinear systems, compartmental systems, network flows, robust control.
\end{keyword}

\end{frontmatter}

	\section{Introduction}\label{sec:introduction}
The study of dynamical flows in infrastructure networks has attracted a considerable amount of attention in recent years. In particular, there is a growing body of literature in the control systems field dealing with issues of stability, optimality, robustness, and resilience in dynamical flow networks. See, e.g., \cite{Paganini:2002,Low.ea:2002,Fan.Arcak.ea:04,ComoPartIITAC13,Bauso.ea:2013,Coogan.Arcak:TAC15,Como:2017} and references therein. 

In this paper, we study deterministic continuous-time models of dynamical flow networks. We consider a finite number of cells exchanging some indistinguishable commodity among themselves and with the external environment. Cells possibly receive a constant exogenous inflow from outside the network and a constant flow is possibly drained out of them directly towards the external environment. We assume that the outflow from a cell is split among its immediately downstream cells in fixed proportions and that each cell has a finite flow and buffer capacity. When the total net flow in a cell ---consisting of the difference between the total flow directed towards it minus the outflow from it--- exceeds the cell's capacity, then the exceeding part of such net flow leaks out of the system. Also, when the difference between the total exogenous demand on a cell and the total inflow in it exceeds the cell's capacity, then the outflow towards the external environment is reduced by an amount equal to the exceeding part of this difference. The ensuing network flow dynamics tuns out to be a linear saturated system with compact state space that we analyse using tools from monotone systems and contraction theory. 

Specifically, we prove that there exists a set of equilibria that is globally asymptotically stable. Such equilibrium set reduces to a single globally asymptotically stable equilibrium for generic exogenous demand vectors. 
Moreover, we show that the critical exogenous demand vectors giving rise to non-unique equilibria correspond to phase transitions in the asymptotic behavior of the dynamical flow network. 

The rest of the paper is organized as follows. The reminder of this section is devoted to the introduction of some notational conventions to be used throughout the paper. In Section \ref{sec:model} we present the class of dynamical flow network models to be studied. Section \ref{sec:main-results} presents the main results concerning the equilibrium set characterization and its global asymptotic stability, as well as the dependence of such equilibria  on the exogenous demand vector. Finally, Section 4 and 5 contain the proofs needed to demonstrate such results.

We shall consider the standard partial order on $\R^n$ whereby the inequality $a\le b$ for two vectors $a,b\in\R^n$ is meant hold true entry-wise. 
A dynamical system with state space $\mc X\subseteq\R^n$ will be referred to as monotone if it preserves such partial order. 
For two vectors $a,b\in\R^n$ such that $a\le b$, we shall denote by 
$$\mathcal{L}_{a}^{b}=\left\{x \in \R^n:\, a \leq x \leq b\right\} =\displaystyle\Pi_{i=1}^n [a_i,b_i]$$
the complete lattice and let $S_a^b:\R^n\to\mc L_a^b$  be the vector saturation function defined by 
\be\label{eq:saturation} \left(S_a^b(y)\right)_i=\max\{a_i,\min\{y_i,b_i\}\}\,,\ee
for $y\in\R^n$ and $i=1,\ldots,n$. For subsets of indices $\mc A,\mc B\subseteq\{1,\ldots,n\}$, we shall denote the restriction of a vector $x\in\R^n$ by $x_{\mc A}=(x_i)_{i\in\mc A}$ and the restriction of a matrix $M\in\R^{n\times n}$ by $M_{\mc A \mc B}=(M_{ij})_{i\in\mc A,j\in\mc B}$. 
 
	\section{A dynamical flow network model with finite capacity}\label{sec:model}
	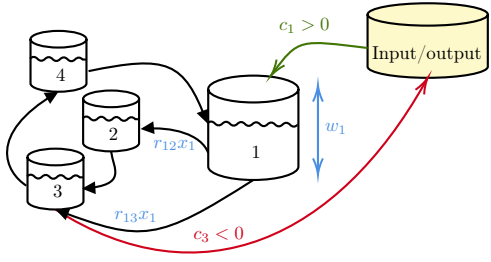
\begin{figure}
	\centering\begin{tikzpicture}[scale=0.75, thick, x=0.75pt,y=0.75pt,yscale=-1,xscale=1]

\draw   (220,131.53) -- (220,183.53) .. controls (220,188.5) and (206.57,192.53) .. (190,192.53) .. controls (173.43,192.53) and (160,188.5) .. (160,183.53) -- (160,131.53)(220,131.53) .. controls (220,136.5) and (206.57,140.53) .. (190,140.53) .. controls (173.43,140.53) and (160,136.5) .. (160,131.53) .. controls (160,126.56) and (173.43,122.53) .. (190,122.53) .. controls (206.57,122.53) and (220,126.56) .. (220,131.53) -- cycle ;
\draw    (159.43,154.53) .. controls (161.13,152.9) and (162.8,152.93) .. (164.43,154.62) .. controls (166.07,156.31) and (167.74,156.34) .. (169.43,154.7) .. controls (171.12,153.06) and (172.79,153.09) .. (174.43,154.78) .. controls (176.07,156.47) and (177.74,156.5) .. (179.43,154.87) .. controls (181.12,153.23) and (182.79,153.26) .. (184.43,154.95) .. controls (186.07,156.64) and (187.74,156.67) .. (189.43,155.03) .. controls (191.12,153.4) and (192.79,153.43) .. (194.43,155.12) .. controls (196.07,156.81) and (197.74,156.84) .. (199.43,155.2) .. controls (201.12,153.56) and (202.79,153.59) .. (204.43,155.28) .. controls (206.07,156.97) and (207.74,157) .. (209.43,155.37) .. controls (211.12,153.73) and (212.79,153.76) .. (214.43,155.45) .. controls (216.06,157.14) and (217.73,157.17) .. (219.42,155.53) -- (219.43,155.53) -- (219.43,155.53) ;

\draw  [fill={rgb, 255:red, 248; green, 231; blue, 28 }  ,fill opacity=0.23 ] (346,81.5) -- (346,116.5) .. controls (346,120.64) and (328.09,124) .. (306,124) .. controls (283.91,124) and (266,120.64) .. (266,116.5) -- (266,81.5)(346,81.5) .. controls (346,85.64) and (328.09,89) .. (306,89) .. controls (283.91,89) and (266,85.64) .. (266,81.5) .. controls (266,77.36) and (283.91,74) .. (306,74) .. controls (328.09,74) and (346,77.36) .. (346,81.5) -- cycle ;
\draw [color={rgb, 255:red, 208; green, 2; blue, 27 }  ,draw opacity=1 ]   (304.25,127.17) .. controls (246.17,231.22) and (155.41,274.7) .. (58.39,212.33) ;

\draw [shift={(306,124)}, rotate = 118.61] [color={rgb, 255:red, 208; green, 2; blue, 27 }  ,draw opacity=1 ][line width=0.75]    (10.93,-3.29) .. controls (6.95,-1.4) and (3.31,-0.3) .. (0,0) .. controls (3.31,0.3) and (6.95,1.4) .. (10.93,3.29)   ;
\draw [color={rgb, 255:red, 65; green, 117; blue, 5 }  ,draw opacity=1 ]   (266,104) .. controls (220.92,98.45) and (220.02,97.38) .. (202.79,121.08) ;
\draw [shift={(201.72,122.57)}, rotate = 305.93] [color={rgb, 255:red, 65; green, 117; blue, 5 }  ,draw opacity=1 ][line width=0.75]    (10.93,-3.29) .. controls (6.95,-1.4) and (3.31,-0.3) .. (0,0) .. controls (3.31,0.3) and (6.95,1.4) .. (10.93,3.29)   ;

\draw   (79,98.12) -- (79,127.75) .. controls (79,130.83) and (70.67,133.33) .. (60.39,133.33) .. controls (50.12,133.33) and (41.78,130.83) .. (41.78,127.75) -- (41.78,98.12)(79,98.12) .. controls (79,101.2) and (70.67,103.7) .. (60.39,103.7) .. controls (50.12,103.7) and (41.78,101.2) .. (41.78,98.12) .. controls (41.78,95.03) and (50.12,92.53) .. (60.39,92.53) .. controls (70.67,92.53) and (79,95.03) .. (79,98.12) -- cycle ;
\draw    (41.43,111.18) .. controls (43.12,109.54) and (44.79,109.57) .. (46.43,111.26) .. controls (48.07,112.95) and (49.74,112.98) .. (51.43,111.34) .. controls (53.12,109.7) and (54.79,109.73) .. (56.43,111.42) .. controls (58.07,113.11) and (59.74,113.14) .. (61.43,111.5) .. controls (63.12,109.86) and (64.79,109.89) .. (66.43,111.58) .. controls (68.08,113.27) and (69.74,113.29) .. (71.43,111.65) .. controls (73.12,110.01) and (74.79,110.04) .. (76.43,111.73) -- (78.65,111.77) -- (78.65,111.77) ;

\draw   (77,177.12) -- (77,206.75) .. controls (77,209.83) and (68.67,212.33) .. (58.39,212.33) .. controls (48.12,212.33) and (39.78,209.83) .. (39.78,206.75) -- (39.78,177.12)(77,177.12) .. controls (77,180.2) and (68.67,182.7) .. (58.39,182.7) .. controls (48.12,182.7) and (39.78,180.2) .. (39.78,177.12) .. controls (39.78,174.03) and (48.12,171.53) .. (58.39,171.53) .. controls (68.67,171.53) and (77,174.03) .. (77,177.12) -- cycle ;
\draw    (39.43,190.18) .. controls (41.12,188.54) and (42.79,188.57) .. (44.43,190.26) .. controls (46.07,191.95) and (47.74,191.98) .. (49.43,190.34) .. controls (51.12,188.7) and (52.79,188.73) .. (54.43,190.42) .. controls (56.07,192.11) and (57.74,192.14) .. (59.43,190.5) .. controls (61.12,188.86) and (62.79,188.89) .. (64.43,190.58) .. controls (66.08,192.27) and (67.74,192.29) .. (69.43,190.65) .. controls (71.12,189.01) and (72.79,189.04) .. (74.43,190.73) -- (76.65,190.77) -- (76.65,190.77) ;

\draw   (114,138.12) -- (114,167.75) .. controls (114,170.83) and (105.67,173.33) .. (95.39,173.33) .. controls (85.12,173.33) and (76.78,170.83) .. (76.78,167.75) -- (76.78,138.12)(114,138.12) .. controls (114,141.2) and (105.67,143.7) .. (95.39,143.7) .. controls (85.12,143.7) and (76.78,141.2) .. (76.78,138.12) .. controls (76.78,135.03) and (85.12,132.53) .. (95.39,132.53) .. controls (105.67,132.53) and (114,135.03) .. (114,138.12) -- cycle ;
\draw    (76.43,151.18) .. controls (78.12,149.54) and (79.79,149.57) .. (81.43,151.26) .. controls (83.07,152.95) and (84.74,152.98) .. (86.43,151.34) .. controls (88.12,149.7) and (89.79,149.73) .. (91.43,151.42) .. controls (93.07,153.11) and (94.74,153.14) .. (96.43,151.5) .. controls (98.12,149.86) and (99.79,149.89) .. (101.43,151.58) .. controls (103.08,153.27) and (104.74,153.29) .. (106.43,151.65) .. controls (108.12,150.01) and (109.79,150.04) .. (111.43,151.73) -- (113.65,151.77) -- (113.65,151.77) ;

\draw    (58.26,134.59) .. controls (12.39,162.11) and (25.21,187.47) .. (39,196.33) ;

\draw [shift={(60.39,133.33)}, rotate = 149.95] [fill={rgb, 255:red, 0; green, 0; blue, 0 }  ][line width=0.75]  [draw opacity=0] (8.93,-4.29) -- (0,0) -- (8.93,4.29) -- cycle    ;
\draw    (95.39,173.33) .. controls (96.97,185.09) and (106.6,195.89) .. (77.81,198.2) ;
\draw [shift={(76,198.33)}, rotate = 356.31] [fill={rgb, 255:red, 0; green, 0; blue, 0 }  ][line width=0.75]  [draw opacity=0] (8.93,-4.29) -- (0,0) -- (8.93,4.29) -- cycle    ;

\draw    (60.95,213.31) .. controls (116.28,234.32) and (125.98,233.71) .. (190,192.53) ;

\draw [shift={(58.39,212.33)}, rotate = 20.9] [fill={rgb, 255:red, 0; green, 0; blue, 0 }  ][line width=0.75]  [draw opacity=0] (8.93,-4.29) -- (0,0) -- (8.93,4.29) -- cycle    ;
\draw    (80.39,119.53) .. controls (133.2,108.5) and (143.69,113.19) .. (158.74,152.7) ;
\draw [shift={(159.43,154.53)}, rotate = 249.45999999999998] [fill={rgb, 255:red, 0; green, 0; blue, 0 }  ][line width=0.75]  [draw opacity=0] (8.93,-4.29) -- (0,0) -- (8.93,4.29) -- cycle    ;

\draw    (160,174.33) .. controls (155.15,161.72) and (136.19,156.64) .. (116.8,158.17) ;
\draw [shift={(115,158.33)}, rotate = 354.28999999999996] [fill={rgb, 255:red, 0; green, 0; blue, 0 }  ][line width=0.75]  [draw opacity=0] (8.93,-4.29) -- (0,0) -- (8.93,4.29) -- cycle    ;

\draw [color={rgb, 255:red, 74; green, 144; blue, 226 }  ,draw opacity=1 ]   (233,131) -- (233,187) ;
\draw [shift={(233,189)}, rotate = 270] [color={rgb, 255:red, 74; green, 144; blue, 226 }  ,draw opacity=1 ][line width=0.75]    (10.93,-3.29) .. controls (6.95,-1.4) and (3.31,-0.3) .. (0,0) .. controls (3.31,0.3) and (6.95,1.4) .. (10.93,3.29)   ;
\draw [shift={(233,129)}, rotate = 90] [color={rgb, 255:red, 74; green, 144; blue, 226 }  ,draw opacity=1 ][line width=0.75]    (10.93,-3.29) .. controls (6.95,-1.4) and (3.31,-0.3) .. (0,0) .. controls (3.31,0.3) and (6.95,1.4) .. (10.93,3.29)   ;

\draw (305.5,108) node [scale=0.7] [align=left] {Input/output};
\draw (245,158) node [scale=0.7,color={rgb, 255:red, 74; green, 144; blue, 226 }  ,opacity=1 ]  {$w_{1}$};
\draw (138,169) node [scale=0.7,color={rgb, 255:red, 74; green, 144; blue, 226 }  ,opacity=1 ]  {$r_{12} x_{1}$};
\draw (192,174) node [scale=0.7]  {$1$};
\draw (97,162) node [scale=0.7]  {$2$};
\draw (60,200) node [scale=0.7]  {$3$};
\draw (62,122) node [scale=0.7]  {$4$};
\draw (112,217) node [scale=0.7,color={rgb, 255:red, 74; green, 144; blue, 226 }  ,opacity=1 ]  {$r_{13} x_{1}$};
\draw (224,89) node [scale=0.7,color={rgb, 255:red, 65; green, 117; blue, 5 }  ,opacity=1 ]  {$c_{1}  >0$};
\draw (167,229) node [scale=0.7,color={rgb, 255:red, 208; green, 2; blue, 27 }  ,opacity=1 ]  {$c_{3} < 0$};

\end{tikzpicture}
\caption{\label{fig:flow-network}Illustration of a dynamical flow network with four cells.}
\end{figure}

We consider dynamical flow networks consisting of finitely many cells $i=1,2,\ldots,n$, exchanging an indistinguishable commodity both among themselves and with the external environment as described below. (See also Figure \ref{fig:flow-network})

Let $x_i(t)$ be the quantity of commodity contained in cell $i=1,2,\ldots,n$ at time $t\ge0$ and let $w_i>0$ be its capacity.
The state of the system is described by the vector $x(t)=(x_i(t))_{1\le i\le n}$ and evolves in continuous time according to the following dynamical system \be\label{flow-dynamics}\dot x=f(x)\,,\ee
where $f(x)=(f_i(x))_{1\le i\le n}$ is the vector of instantaneous net flows (inflows minus outflows) in the cells that will be assumed to satisfy the constraints 
\be\label{flow-constraint}-x_i\le f_i(x)\le w_i-x_i\,,\qquad i=0,\ldots,n\,,\ee
throughout the evolution of the system. 
Notice that the leftmost inequality in \eqref{flow-constraint} states that the outflow from cell $i$ can never exceed the current inflow plus the total quantity of commodity in the cell, in particular implying the physically meaningful fact that the net flow $f_i(x)$ is nonnegative when the cell is empty (i.e., when $x_i=0$) so that $x_i(t)$ can never become negative. On the other hand, the rightmost inequality in \eqref{flow-constraint} guarantees that the sum of the current total mass and the inflow in a cell $i$ and can never exceed the difference between its capacity $w_i$ and the current outflow, so that in particular, when the mass $x(t)=w_i$ has reached the capacity, the net flow $f_i(x)$ is nonpositive, thus implying that the total mass will never exceed the capacity $w_i$ if started below that. The complete lattice $\mathcal{L}_{0}^{w}$ is invariant for any dynamical flow network \eqref{flow-dynamics} satisfying \eqref{flow-constraint}. 

Now, let each cell $i$ possibly receive a constant exogenous inflow $\lambda_i\ge0$ from outside the network and let a constant flow $\mu_i\ge0$ possibly be drained directly from cell $i$ towards the external environment, 
and let $c_i=\lambda_i-\mu_i$ be the exogenous net demand on cell $i$. 
Also, assume that constant fraction $R_{ij}\ge0$ of the quantity of commodity $x_i$ flows directly towards another cell $j\ne i$ in the network (fixed routing), while the remaining part $(1-\sum_jR_{ij})x_i$ leaves the network directly. Notice that the routing matrix $R=(R_{ij})\in\R^{n\times n}$  is necessarily sub-stochastic, i.e., with nonnegative entries and such that its rows all have sum less than or equal to $1$. 

Conservation of mass and the constraint \eqref{flow-constraint} imply that the netflow in each cell $i=1,\ldots,n$ is given by 
\be\label{flow-exact}\ba{rcl} 
f_i(x)&=& S_{-x_i}^{w_i-x_i}\left(\lambda_i-\mu_i+\sum\nolimits_jR_{ji}x_j-x_i\right)\\[3pt]
&=&S_0^{w_i}\left(\sum\nolimits_jR_{ji}x_j+c_i\right)-x_i\,.\ea\ee 
We may then rewrite the dynamical flow network \eqref{flow-dynamics}--\eqref{flow-exact} compactly as 
\begin{equation}\label{ema}
\dot x  = S_0^{w} \left( R'x+c\right)-x\,,
\end{equation}
where $w\in\R^n$ is the vector of the cells' capacities. Observe that the function 
$f(x)$ as defined in \eqref{flow-exact} is Lipschitz continuous in $\R^n$, so that existence and uniqueness of a solution to the dynamical flow network \eqref{ema} is ensured for every initial state $x(0)\in\mc L_0^w$. 

Observe that in the dynamical network flow \eqref{ema} it is understood that when the difference between the total flow $\lambda_i+\sum_jR_{ji}x_j$ directed towards a cell and the outflow $\mu_i+x_i$ from it exceeds the capacity $w_i$, then the exceeding part of it leaks out of the system. Moreover, the dynamical network flow \eqref{ema} also assumes that, when the difference between the total exogenous demand $\mu_i$ on a cell $i$ and the total inflow $\lambda_i+\sum_jR_{ji}x_j$ exceeds the cell's capacity $w_i$, then the outflow towards the external environment is reduced by an amount equal to the exceeding part of this difference.


\section{Main results}\label{sec:main-results}  
In this section, we state the main results of this paper. These are concerned on the one hand with global asymptotic stability of the dynamical flow network \eqref{ema} and on the other hand on the dependance (in particular, continuity and the lack thereof) of the equilibria of \eqref{ema} on the exogenous demand vector $c\in\R^n$. 

Before proceeding, let us gather some terminology that is used in our statements. The routing matrix $R$ will be referred to as out-connected is for every $i=1,\ldots,n$ there exists $j\in\{1,\ldots,n\}$ such that $\sum_{k}R_{jk}<1$ and $(R^l)_{ij}>0$ for some $l\ge0$. It will be referred to as stochastic if all its rows sum up to $1$ and irreducible if, for every nonempty proper subset $\mc S\subsetneq\{1,\ldots,n\}$, there exists at least one $i\in\mc S$ such that $\sum_{j\in\mc S}R_{ij}<1$. 
It is a standard fact that, if the routing matrix $R$ is stochastic irreducible, then it admits a unique invariant probability vector $\pi=R'\pi$ and such vector is strictly positive entry-wise. Moreover, for every zero-sum vector $v\in\R^n$, the vector series 
\begin{equation}\label{eq:series}
Hv :=  \frac{1}{2} \sum_{k \geq 0}\left(\frac{I+R^{\prime}}{2}\right)^{k} v
\end{equation}
is convergent and its limit satisfies 
\be\label{Hv=R'Hv}Hv=R'Hv+v\,.\ee

We start with the stability results that are stated in the following. 
\begin{theorem}\label{TH:UNIQ}
Let $w\in\R^n$ be a positive vector and $R\in\R^{n\times n}$ a sub-stochastic matrix. Then, 
\begin{enumerate}
\item[(i)] if $R$ is sub-stochastic and out-connected, then, for every exogenous demand vector $c\in\R^n$ the dynamical flow network \eqref{ema} admits a globally asymptotically stable equilibrium $x^*\in\mc L^w_0$. 
\end{enumerate}
On the other hand, if $R$ is stochastic and irreducible, then
\begin{enumerate}
\item[(ii)]  for every exogenous demand vector $c\in\R^n$ the set of equilibria $\mc X(c)$ of the dynamical flow network \eqref{ema} is a nonempty line segment joining two points $\ul x\le\ov x$ on the boundary of the lattice $\mc L_0^w$;  
\item[(iii)] for every initial state $x(0)\in\mc L_0^w$, the solution of  \eqref{ema} converges to the equilibrium set $\mc X(c)$ as $t$ grows large; 
\item[(iv)]  the equilibrium set $\mc X(c)$ has positive length if and only if 
\begin{equation}\label{eq:nonuniqe}
\min _{i}\left\{\frac{(Hc)_{i}}{\pi_{i}}\right\}+\min _{i}\left\{\frac{w_{i}-(Hc)_{i}}{\pi_{i}}\right\}>0
\end{equation}
\end{enumerate}
\end{theorem}


Theorem 1 characterizes the set of equilibria $\mc X(c)$ and it is particularly relevant, in a given network, to study the behavior of such set with respect to possible variations of the exogenous net flow vector $c$. Indeed, this exogenous flow might be subject to shocks and variations that might affect the whole flow on the network. Thus, the resilience of the system with respect to shocks is in the end determined by the way solutions depend on the parameter vector $c$. We will show that there exists a set of critical vector $c$ such that the equilibria of \eqref{ema} undergo a jump discontinuity, thus determining a phase transition in the asymptotic behavior of the system, and we will describe this critical set.

 Let us introduce some notation that will be used in the main statement. 
Let 
\be\label{eq:UM-def}\mc U=\left\{c\in \R^n:\,|\mc X (c)|=1\right\}\,,\quad \mc M=\R^n\setminus\mc  U\,,\ee
be the subsets of exogenous flow vectors for which there is a unique equilibrium and, respectively, there are multiple equilibria. Moreover, we denote with $\ul x(c)$ and $\ov x(c)$ the smallest and largest equilibria for a given vector $c$.
For exogenous flow vectors $c\in\mc U$, we shall also use the notation $$x^*(c)=\ul x(c)=\ov x(c)$$ 
for the unique equilibrium. 

We can now state the following result.

\begin{theorem}\label{TH:CONTINUITY}
	Let $w\in\R_+^n$ be a nonnegative vector. 
	Let $\mc U$ and $\mc M$ be defined as in \eqref{eq:UM-def}.
	Then, 
	\begin{enumerate}
		\item[(i)] if $R$ is sub-stochastic and out-connected, then, for every exogenous demand vector $c\in\R^n$ the map $c\mapsto x^*(c)$ is continuous.
	\end{enumerate}
On the other hand, if $R$ is stochastic and irreducible, then
	\begin{enumerate}
		\item[(ii)] $\mc  M$ is linear sub-manifold of co-dimension $1$;
		\item[(iii)] the map $c\mapsto x^*(c)$ is continuous on the set $\mc U$; 
		\item[(iv)] for every $c^*\in\mc M$, 
		$$ \liminf\limits_{{\substack{c\in\mc U\\ c\to c^*}}} x^*(c)=\underline x(c^*)\,,\qquad \limsup\limits_{\substack{c\in\mc U\\ c\to c^*}}x^*(c)=\bar x(c^*)\,.$$
	\end{enumerate}
\end{theorem}
Theorem \ref{TH:CONTINUITY}, and in particular the condition (iv), states that the equilibria of $\eqref{ema}$ undergo a jump discontinuity when the vector $c$ crosses the set $\mc M$ for which the uniqueness condition for equilibria fails to hold. This in turn implies that even a slight change in the exogenous flow may trigger a phase transition in the system and a huge impact on the quantity of commodities exchanged at equilibrium in the network. We show this phenomenon in the following example.

\subsubsection{Example 1}
\begin{figure}[H]
	\centering	
	\begin{tikzpicture}[scale=0.85, x=0.75pt,y=0.75pt,yscale=-1,xscale=1]
	
	\draw   (108.57,174.28) .. controls (108.57,165.6) and (115.6,158.57) .. (124.28,158.57) .. controls (132.96,158.57) and (140,165.6) .. (140,174.28) .. controls (140,182.96) and (132.96,190) .. (124.28,190) .. controls (115.6,190) and (108.57,182.96) .. (108.57,174.28) -- cycle ;
	\draw   (160,114.28) .. controls (160,105.6) and (167.04,98.57) .. (175.72,98.57) .. controls (184.4,98.57) and (191.43,105.6) .. (191.43,114.28) .. controls (191.43,122.96) and (184.4,130) .. (175.72,130) .. controls (167.04,130) and (160,122.96) .. (160,114.28) -- cycle ;
	\draw   (210,174.28) .. controls (210,165.6) and (217.04,158.57) .. (225.72,158.57) .. controls (234.4,158.57) and (241.43,165.6) .. (241.43,174.28) .. controls (241.43,182.96) and (234.4,190) .. (225.72,190) .. controls (217.04,190) and (210,182.96) .. (210,174.28) -- cycle ;
	\draw    (124.28,158.57) .. controls (119.55,131.3) and (133.7,115.41) .. (158.1,114.34) ;
	\draw [shift={(160,114.28)}, rotate = 539.29] [fill={rgb, 255:red, 0; green, 0; blue, 0 }  ][line width=0.75]  [draw opacity=0] (8.93,-4.29) -- (0,0) -- (8.93,4.29) -- cycle    ;
	
	\draw    (225.72,158.57) .. controls (235.24,127.24) and (223.9,114.13) .. (193.33,114.26) ;
	\draw [shift={(191.43,114.28)}, rotate = 358.78] [fill={rgb, 255:red, 0; green, 0; blue, 0 }  ][line width=0.75]  [draw opacity=0] (8.93,-4.29) -- (0,0) -- (8.93,4.29) -- cycle    ;
	
	\draw    (175.72,130) .. controls (179.38,157.19) and (181.37,155.65) .. (208.73,173.46) ;
	\draw [shift={(210,174.28)}, rotate = 213.19] [fill={rgb, 255:red, 0; green, 0; blue, 0 }  ][line width=0.75]  [draw opacity=0] (8.93,-4.29) -- (0,0) -- (8.93,4.29) -- cycle    ;
	
	\draw    (175.72,130) .. controls (171.5,155.22) and (171.43,156.56) .. (141.39,173.5) ;
	\draw [shift={(140,174.28)}, rotate = 330.64] [fill={rgb, 255:red, 0; green, 0; blue, 0 }  ][line width=0.75]  [draw opacity=0] (8.93,-4.29) -- (0,0) -- (8.93,4.29) -- cycle    ;
	
	\draw    (124.28,190) .. controls (173.93,203.46) and (175.41,204.58) .. (224.22,190.43) ;
	\draw [shift={(225.72,190)}, rotate = 523.81] [fill={rgb, 255:red, 0; green, 0; blue, 0 }  ][line width=0.75]  [draw opacity=0] (8.93,-4.29) -- (0,0) -- (8.93,4.29) -- cycle    ;
	
	\draw  [fill={rgb, 255:red, 248; green, 231; blue, 28 }  ,fill opacity=0.23 ] (320,57.5) -- (320,92.5) .. controls (320,96.64) and (302.09,100) .. (280,100) .. controls (257.91,100) and (240,96.64) .. (240,92.5) -- (240,57.5)(320,57.5) .. controls (320,61.64) and (302.09,65) .. (280,65) .. controls (257.91,65) and (240,61.64) .. (240,57.5) .. controls (240,53.36) and (257.91,50) .. (280,50) .. controls (302.09,50) and (320,53.36) .. (320,57.5) -- cycle ;
	\draw [color={rgb, 255:red, 208; green, 2; blue, 27 }  ,draw opacity=1 ]   (280.5,102.57) .. controls (288.45,143.65) and (288.47,145.19) .. (241.43,174.28) ;
	
	\draw [shift={(280,100)}, rotate = 79.05] [color={rgb, 255:red, 208; green, 2; blue, 27 }  ,draw opacity=1 ][line width=0.75]    (10.93,-3.29) .. controls (6.95,-1.4) and (3.31,-0.3) .. (0,0) .. controls (3.31,0.3) and (6.95,1.4) .. (10.93,3.29)   ;
	\draw [color={rgb, 255:red, 65; green, 117; blue, 5 }  ,draw opacity=1 ]   (240,80) .. controls (202.01,73.7) and (201.44,73.6) .. (176.86,97.46) ;
	\draw [shift={(175.72,98.57)}, rotate = 315.85] [color={rgb, 255:red, 65; green, 117; blue, 5 }  ,draw opacity=1 ][line width=0.75]    (10.93,-3.29) .. controls (6.95,-1.4) and (3.31,-0.3) .. (0,0) .. controls (3.31,0.3) and (6.95,1.4) .. (10.93,3.29)   ;

	\draw (124.28,174.28) node   {$1$};
	\draw (175.72,114.28) node   {$3$};
	\draw (225.72,174.28) node   {$2$};
	\draw (107.5,127) node [scale=0.7]  {$0.25$};
	\draw (237,123) node [scale=0.7]  {$1$};
	\draw (197.5,150) node [scale=0.7]  {$0.7$};
	\draw (150.5,153) node [scale=0.7]  {$0.3$};
	\draw (172.5,193) node [scale=0.7]  {$0.75$};
	\draw (279.5,84) node [scale=0.7] [align=left] {Input/output};
	\draw (298.5,137) node [scale=0.7,color={rgb, 255:red, 208; green, 2; blue, 27 }  ,opacity=1 ]  {$-1$};
	\draw (203,67) node [scale=0.7,color={rgb, 255:red, 65; green, 117; blue, 5 }  ,opacity=1 ]  {$1$};

	\end{tikzpicture}
	
	\caption{Flow network with three cells.}
	\label{f0N2}
\end{figure}
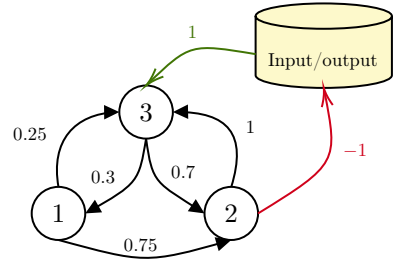

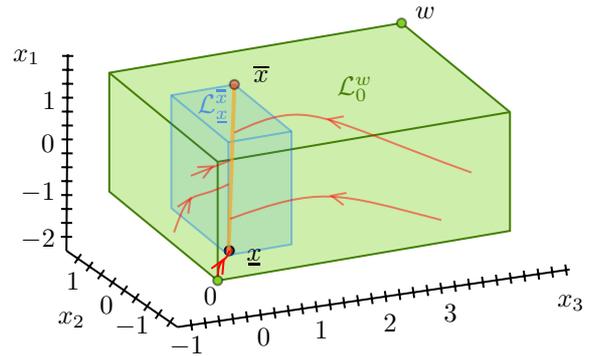
\begin{figure}[H]
	\centering
	\begin{tikzpicture}[scale=0.7, x=0.75pt,y=0.75pt,yscale=-1,xscale=1]
	
	\draw  [color={rgb, 255:red, 69; green, 127; blue, 6 }  ,draw opacity=1 ][fill={rgb, 255:red, 126; green, 211; blue, 33 }  ,fill opacity=0.38 ] (252.13,212.68) -- (174.84,148.66) -- (174.83,63.31) -- (383.16,27.55) -- (460.45,91.56) -- (460.46,176.92) -- cycle ; \draw  [color={rgb, 255:red, 69; green, 127; blue, 6 }  ,draw opacity=1 ] (174.83,63.31) -- (252.12,127.32) -- (252.13,212.68) ; \draw  [color={rgb, 255:red, 69; green, 127; blue, 6 }  ,draw opacity=1 ] (252.12,127.32) -- (460.45,91.56) ;
	\draw    (144.3,191.2) -- (145.3,50.2) (140.37,181.17) -- (148.37,181.23)(140.44,171.17) -- (148.44,171.23)(140.51,161.17) -- (148.51,161.23)(140.58,151.17) -- (148.58,151.23)(140.65,141.17) -- (148.65,141.23)(140.73,131.17) -- (148.73,131.23)(140.8,121.17) -- (148.8,121.23)(140.87,111.17) -- (148.87,111.23)(140.94,101.17) -- (148.94,101.23)(141.01,91.17) -- (149.01,91.23)(141.08,81.17) -- (149.08,81.23)(141.15,71.17) -- (149.15,71.23)(141.22,61.17) -- (149.22,61.23)(141.29,51.18) -- (149.29,51.23) ;

	\draw    (144.3,191.2) -- (223.47,246.06) (154.8,193.61) -- (150.24,200.18)(163.02,199.3) -- (158.46,205.88)(171.24,205) -- (166.68,211.58)(179.46,210.7) -- (174.9,217.27)(187.68,216.39) -- (183.12,222.97)(195.89,222.09) -- (191.34,228.66)(204.11,227.78) -- (199.56,234.36)(212.33,233.48) -- (207.78,240.05)(220.55,239.17) -- (216,245.75) ;

	\draw    (223.47,246.06) -- (503,204.67) (232.77,240.64) -- (233.94,248.55)(242.66,239.17) -- (243.84,247.09)(252.56,237.71) -- (253.73,245.62)(262.45,236.24) -- (263.62,244.16)(272.34,234.78) -- (273.51,242.69)(282.23,233.31) -- (283.4,241.23)(292.12,231.85) -- (293.3,239.76)(302.02,230.38) -- (303.19,238.3)(311.91,228.92) -- (313.08,236.83)(321.8,227.45) -- (322.97,235.37)(331.69,225.99) -- (332.86,233.9)(341.59,224.53) -- (342.76,232.44)(351.48,223.06) -- (352.65,230.97)(361.37,221.6) -- (362.54,229.51)(371.26,220.13) -- (372.43,228.04)(381.15,218.67) -- (382.33,226.58)(391.05,217.2) -- (392.22,225.11)(400.94,215.74) -- (402.11,223.65)(410.83,214.27) -- (412,222.19)(420.72,212.81) -- (421.89,220.72)(430.61,211.34) -- (431.79,219.26)(440.51,209.88) -- (441.68,217.79)(450.4,208.41) -- (451.57,216.33)(460.29,206.95) -- (461.46,214.86)(470.18,205.48) -- (471.35,213.4)(480.08,204.02) -- (481.25,211.93)(489.97,202.55) -- (491.14,210.47)(499.86,201.09) -- (501.03,209) ;

	\draw  [color={rgb, 255:red, 0; green, 0; blue, 0 }  ,draw opacity=1 ][fill={rgb, 255:red, 255; green, 0; blue, 0 }  ,fill opacity=1 ] (257.21,190.55) .. controls (257.54,188.79) and (259.25,187.63) .. (261.01,187.97) .. controls (262.77,188.31) and (263.93,190.01) .. (263.59,191.77) .. controls (263.25,193.54) and (261.55,194.69) .. (259.79,194.36) .. controls (258.02,194.02) and (256.87,192.32) .. (257.21,190.55) -- cycle ;
	\draw  [color={rgb, 255:red, 0; green, 118; blue, 255 }  ,draw opacity=0.39 ][fill={rgb, 255:red, 33; green, 146; blue, 211 }  ,fill opacity=0.14 ] (259.79,194.36) -- (219.06,160.62) -- (219.05,79.45) -- (264.02,71.73) -- (304.75,105.47) -- (304.76,186.64) -- cycle ; \draw  [color={rgb, 255:red, 0; green, 118; blue, 255 }  ,draw opacity=0.39 ] (219.05,79.45) -- (259.77,113.19) -- (259.79,194.36) ; \draw  [color={rgb, 255:red, 0; green, 118; blue, 255 }  ,draw opacity=0.39 ] (259.77,113.19) -- (304.75,105.47) ;
	\draw [color={rgb, 255:red, 245; green, 166; blue, 35 }  ,draw opacity=0.69 ][line width=1.5]    (264.02,71.73) -- (259.79,194.36) ;

	\draw  [color={rgb, 255:red, 0; green, 0; blue, 0 }  ,draw opacity=0.44 ][fill={rgb, 255:red, 255; green, 0; blue, 0 }  ,fill opacity=0.39 ] (260.83,71.12) .. controls (261.17,69.36) and (262.87,68.2) .. (264.63,68.54) .. controls (266.39,68.88) and (267.55,70.58) .. (267.21,72.34) .. controls (266.88,74.11) and (265.17,75.26) .. (263.41,74.93) .. controls (261.65,74.59) and (260.49,72.89) .. (260.83,71.12) -- cycle ;
	\draw [color={rgb, 255:red, 255; green, 0; blue, 0 }  ,draw opacity=1 ][line width=0.75]    (252.13,212.68) .. controls (252,194.33) and (258,201.33) .. (260.4,191.16) ;

	\draw [color={rgb, 255:red, 255; green, 0; blue, 0 }  ,draw opacity=0.48 ][line width=0.75]    (233,137.33) -- (261,126.33) ;

	\draw [color={rgb, 255:red, 255; green, 0; blue, 0 }  ,draw opacity=0.48 ]   (221,178) .. controls (233,143.33) and (232,157.33) .. (260,143.33) ;

	\draw [color={rgb, 255:red, 255; green, 0; blue, 0 }  ,draw opacity=0.48 ]   (264,106.33) .. controls (318,80.33) and (323,95.33) .. (433,135) ;

	\draw [color={rgb, 255:red, 255; green, 0; blue, 0 }  ,draw opacity=0.48 ]   (261,168.33) .. controls (326,144.33) and (336,149.33) .. (411.5,169.17) ;

	\draw  [color={rgb, 255:red, 255; green, 0; blue, 0 }  ,draw opacity=0.48 ] (341.82,158.37) -- (333.14,151.67) -- (343.91,149.62) ;
	\draw  [color={rgb, 255:red, 255; green, 0; blue, 0 }  ,draw opacity=0.48 ] (339.82,103.37) -- (331.14,96.67) -- (341.91,94.62) ;
	\draw  [color={rgb, 255:red, 255; green, 0; blue, 0 }  ,draw opacity=0.48 ] (237.49,131.14) -- (248.44,130.53) -- (241.63,139.13) ;
	\draw  [color={rgb, 255:red, 255; green, 0; blue, 0 }  ,draw opacity=0.48 ] (223.52,159.64) -- (232.68,153.61) -- (231.12,164.47) ;
	\draw  [color={rgb, 255:red, 255; green, 0; blue, 0 }  ,draw opacity=1 ] (247.63,203.99) -- (256.45,197.47) -- (255.47,208.4) ;
	\draw  [color={rgb, 255:red, 0; green, 0; blue, 0 }  ,draw opacity=0.44 ][fill={rgb, 255:red, 126; green, 211; blue, 33 }  ,fill opacity=1 ] (379.97,26.94) .. controls (380.3,25.18) and (382.01,24.02) .. (383.77,24.36) .. controls (385.53,24.69) and (386.69,26.4) .. (386.35,28.16) .. controls (386.01,29.92) and (384.31,31.08) .. (382.55,30.74) .. controls (380.79,30.4) and (379.63,28.7) .. (379.97,26.94) -- cycle ;
	\draw  [color={rgb, 255:red, 0; green, 0; blue, 0 }  ,draw opacity=0.44 ][fill={rgb, 255:red, 126; green, 211; blue, 33 }  ,fill opacity=1 ] (248.94,212.07) .. controls (249.28,210.3) and (250.98,209.15) .. (252.74,209.48) .. controls (254.51,209.82) and (255.66,211.52) .. (255.33,213.29) .. controls (254.99,215.05) and (253.29,216.21) .. (251.52,215.87) .. controls (249.76,215.53) and (248.6,213.83) .. (248.94,212.07) -- cycle ;
	
	\draw (148,241) node   {$x_{2}$};
	\draw (504,229) node   {$x_{3}$};
	\draw (131,114) node   {$0$};
	\draw (124,149) node   {$-1$};
	\draw (124,183) node   {$-2$};
	\draw (326,248) node   {$1$};
	\draw (285,253) node   {$0$};
	\draw (230,260) node   {$-1$};
	\draw (349.89,74.82) node [color={rgb, 255:red, 69; green, 125; blue, 3 }  ,opacity=1 ,rotate=-358.52]  {$\mathcal{L}_0^w$};
	\draw (278.17,196.2) node [rotate=-358.52]  {$\underbarr{x}$};
	\draw (248.89,87.82) node [color={rgb, 255:red, 74; green, 144; blue, 226 }  ,opacity=1 ,rotate=-358.52]  {$\mathcal{L}_{\underbarr{x}}^{\overbar{x}}$};
	\draw (283.17,64.2) node [rotate=-358.52]  {$\overline{x}$};
	\draw (375,243) node   {$2$};
	\draw (418,236) node   {$3$};
	\draw (173,230) node   {$0$};
	\draw (149,213) node   {$1$};
	\draw (191,245) node   {$-1$};
	\draw (132,83) node   {$1$};
	\draw (116,52) node   {$x_{1}$};
	\draw (400,20) node   {$w$};
	\draw (247,224) node   {$0$};

	\end{tikzpicture}
	
	\caption{Trajectories in the phase space in case of multiple equilibria.}
	\label{sp1}

\end{figure}

Let us consider a flow model with an irreducible routing matrix $R$, in particular, we consider \eqref{ema} with:
\begin{equation*}\label{key}R=
\begin{bmatrix}
0 &0.75  &0.25 \\ 
0&  0& 1 \\ 
0.3& 0.7 & 0
\end{bmatrix}, \: w=\begin{bmatrix}
5 \\ 4\\ 6
\end{bmatrix},
\:  c=\begin{bmatrix}
0 \\-1 \\ 1
\end{bmatrix}
\end{equation*}
The corresponding flow network is shown in Fig. \ref{f0N2}. 
Since $\1'c=0$ and $	\min _{i}\left\{\frac{(Hc)_{i}}{\pi_{i}}\right\}+\min _{i}\left\{\frac{w_{i}-(Hc)_{i}}{\pi_{i}}\right\} \approx 9.62>0$ then \eqref{ema} admits multiple equilibria because of Theorem \ref{TH:UNIQ}(iv). Indeed one can compute $\overbar{x}(c) \approx [1.62, 4 , 5.41]'$ and $\underbarr{x}(c) \approx [0.32, 0, 1.08]'$. We highlight the big jump that occurs for this particular vector $c$; notice how in the largest solution $\overbar{x}$, cell 2 can deliver its total outflow capacity $4$ while in the smallest solution $\underbarr{x}$ it outputs $0$. A slight change of the exogenous flow around $c$ could then have a huge impact on the network.  In Fig. \ref{sp1} we show some trajectories (in red) for different initial conditions  in the phase space; we also plot the two lattices $\mc L_0^w$ and $\mc L_{\underbarr{x}}^{\overbar{x}}$ (in green and light blue respectively); finally, the segment of equilibria $ \mc X $ is plot in orange.
We can notice how all trajectories (red curves) converge to the set of equilibria (orange segment).

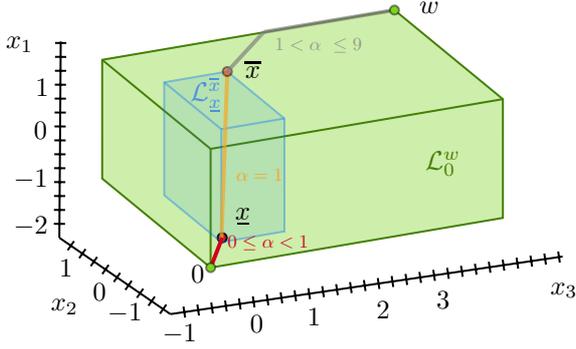
\begin{figure}[H]
	\centering
	\begin{tikzpicture}[scale=0.7 ,x=0.75pt,y=0.75pt,yscale=-1,xscale=1]
	
	\draw  [color={rgb, 255:red, 69; green, 127; blue, 6 }  ,draw opacity=1 ][fill={rgb, 255:red, 126; green, 211; blue, 33 }  ,fill opacity=0.38 ] (246.13,200.68) -- (168.84,136.66) -- (168.83,51.31) -- (377.16,15.55) -- (454.45,79.56) -- (454.46,164.92) -- cycle ; \draw  [color={rgb, 255:red, 69; green, 127; blue, 6 }  ,draw opacity=1 ] (168.83,51.31) -- (246.12,115.32) -- (246.13,200.68) ; \draw  [color={rgb, 255:red, 69; green, 127; blue, 6 }  ,draw opacity=1 ] (246.12,115.32) -- (454.45,79.56) ;
	\draw    (138.3,179.2) -- (139.3,38.2) (134.37,169.17) -- (142.37,169.23)(134.44,159.17) -- (142.44,159.23)(134.51,149.17) -- (142.51,149.23)(134.58,139.17) -- (142.58,139.23)(134.65,129.17) -- (142.65,129.23)(134.73,119.17) -- (142.73,119.23)(134.8,109.17) -- (142.8,109.23)(134.87,99.17) -- (142.87,99.23)(134.94,89.17) -- (142.94,89.23)(135.01,79.17) -- (143.01,79.23)(135.08,69.17) -- (143.08,69.23)(135.15,59.17) -- (143.15,59.23)(135.22,49.17) -- (143.22,49.23)(135.29,39.18) -- (143.29,39.23) ;

	\draw    (138.3,179.2) -- (217.47,234.06) (148.8,181.61) -- (144.24,188.18)(157.02,187.3) -- (152.46,193.88)(165.24,193) -- (160.68,199.58)(173.46,198.7) -- (168.9,205.27)(181.68,204.39) -- (177.12,210.97)(189.89,210.09) -- (185.34,216.66)(198.11,215.78) -- (193.56,222.36)(206.33,221.48) -- (201.78,228.05)(214.55,227.17) -- (210,233.75) ;

	\draw    (217.47,234.06) -- (497,192.67) (226.77,228.64) -- (227.94,236.55)(236.66,227.17) -- (237.84,235.09)(246.56,225.71) -- (247.73,233.62)(256.45,224.24) -- (257.62,232.16)(266.34,222.78) -- (267.51,230.69)(276.23,221.31) -- (277.4,229.23)(286.12,219.85) -- (287.3,227.76)(296.02,218.38) -- (297.19,226.3)(305.91,216.92) -- (307.08,224.83)(315.8,215.45) -- (316.97,223.37)(325.69,213.99) -- (326.86,221.9)(335.59,212.53) -- (336.76,220.44)(345.48,211.06) -- (346.65,218.97)(355.37,209.6) -- (356.54,217.51)(365.26,208.13) -- (366.43,216.04)(375.15,206.67) -- (376.33,214.58)(385.05,205.2) -- (386.22,213.11)(394.94,203.74) -- (396.11,211.65)(404.83,202.27) -- (406,210.19)(414.72,200.81) -- (415.89,208.72)(424.61,199.34) -- (425.79,207.26)(434.51,197.88) -- (435.68,205.79)(444.4,196.41) -- (445.57,204.33)(454.29,194.95) -- (455.46,202.86)(464.18,193.48) -- (465.35,201.4)(474.08,192.02) -- (475.25,199.93)(483.97,190.55) -- (485.14,198.47)(493.86,189.09) -- (495.03,197) ;

	\draw  [color={rgb, 255:red, 0; green, 0; blue, 0 }  ,draw opacity=1 ][fill={rgb, 255:red, 255; green, 0; blue, 0 }  ,fill opacity=1 ] (251.21,178.55) .. controls (251.54,176.79) and (253.25,175.63) .. (255.01,175.97) .. controls (256.77,176.31) and (257.93,178.01) .. (257.59,179.77) .. controls (257.25,181.54) and (255.55,182.69) .. (253.79,182.36) .. controls (252.02,182.02) and (250.87,180.32) .. (251.21,178.55) -- cycle ;
	\draw  [color={rgb, 255:red, 0; green, 118; blue, 255 }  ,draw opacity=0.39 ][fill={rgb, 255:red, 33; green, 146; blue, 211 }  ,fill opacity=0.14 ] (253.79,182.36) -- (213.06,148.62) -- (213.05,67.45) -- (258.02,59.73) -- (298.75,93.47) -- (298.76,174.64) -- cycle ; \draw  [color={rgb, 255:red, 0; green, 118; blue, 255 }  ,draw opacity=0.39 ] (213.05,67.45) -- (253.77,101.19) -- (253.79,182.36) ; \draw  [color={rgb, 255:red, 0; green, 118; blue, 255 }  ,draw opacity=0.39 ] (253.77,101.19) -- (298.75,93.47) ;
	\draw [color={rgb, 255:red, 245; green, 166; blue, 35 }  ,draw opacity=0.69 ][line width=1.5]    (258.02,59.73) -- (253.79,182.36) ;

	\draw  [color={rgb, 255:red, 0; green, 0; blue, 0 }  ,draw opacity=0.44 ][fill={rgb, 255:red, 255; green, 0; blue, 0 }  ,fill opacity=0.39 ] (254.83,59.12) .. controls (255.17,57.36) and (256.87,56.2) .. (258.63,56.54) .. controls (260.39,56.88) and (261.55,58.58) .. (261.21,60.34) .. controls (260.88,62.11) and (259.17,63.26) .. (257.41,62.93) .. controls (255.65,62.59) and (254.49,60.89) .. (254.83,59.12) -- cycle ;
	\draw [color={rgb, 255:red, 208; green, 2; blue, 27 }  ,draw opacity=1 ][line width=1.5]    (254.4,179.16) -- (246.13,200.68) ;

	\draw [color={rgb, 255:red, 128; green, 128; blue, 128 }  ,draw opacity=0.69 ][line width=1.5]    (284,31.33) -- (258.02,59.73) ;

	\draw [color={rgb, 255:red, 128; green, 128; blue, 128 }  ,draw opacity=0.69 ][line width=1.5]    (284,31.33) -- (377.16,15.55) ;

	\draw  [color={rgb, 255:red, 0; green, 0; blue, 0 }  ,draw opacity=0.44 ][fill={rgb, 255:red, 126; green, 211; blue, 33 }  ,fill opacity=1 ] (373.97,14.94) .. controls (374.3,13.18) and (376.01,12.02) .. (377.77,12.36) .. controls (379.53,12.69) and (380.69,14.4) .. (380.35,16.16) .. controls (380.01,17.92) and (378.31,19.08) .. (376.55,18.74) .. controls (374.79,18.4) and (373.63,16.7) .. (373.97,14.94) -- cycle ;
	\draw  [color={rgb, 255:red, 0; green, 0; blue, 0 }  ,draw opacity=0.44 ][fill={rgb, 255:red, 126; green, 211; blue, 33 }  ,fill opacity=1 ] (242.94,200.07) .. controls (243.28,198.3) and (244.98,197.15) .. (246.74,197.48) .. controls (248.51,197.82) and (249.66,199.52) .. (249.33,201.29) .. controls (248.99,203.05) and (247.29,204.21) .. (245.52,203.87) .. controls (243.76,203.53) and (242.6,201.83) .. (242.94,200.07) -- cycle ;
	
	\draw (142,229) node   {$x_{2}$};
	\draw (498,217) node   {$x_{3}$};
	\draw (125,102) node   {$0$};
	\draw (118,137) node   {$-1$};
	\draw (118,171) node   {$-2$};
	\draw (320,236) node   {$1$};
	\draw (279,241) node   {$0$};
	\draw (224,248) node   {$-1$};
	\draw (411.89,125.82) node [color={rgb, 255:red, 69; green, 125; blue, 3 }  ,opacity=1 ,rotate=-358.52]  {$\mathcal{L}_0^w$};
	\draw (269.17,164.2) node [rotate=-358.52]  {$\underbarr{x}$};
	\draw (242.89,76.82) node [color={rgb, 255:red, 74; green, 144; blue, 226 }  ,opacity=1 ,rotate=-358.52]  {$\mathcal{L}_{\underbarr{x}}^{\overbar{x}}$};
	\draw (276.17,58.2) node [rotate=-358.52]  {$\overline{x}$};
	\draw (369,231) node   {$2$};
	\draw (412,224) node   {$3$};
	\draw (167,218) node   {$0$};
	\draw (143,201) node   {$1$};
	\draw (185,233) node   {$-1$};
	\draw (126,71) node   {$1$};
	\draw (110,40) node   {$x_{1}$};
	\draw (402,14) node   {$w$};
	\draw (237,205) node   {$0$};
	\draw (281,134) node [scale=0.7,color={rgb, 255:red, 245; green, 166; blue, 35 }  ,opacity=1 ]  {$\alpha =1$};
	\draw (287,183) node [scale=0.7,color={rgb, 255:red, 208; green, 2; blue, 27 }  ,opacity=1 ]  {$0\leq \alpha < 1$};
	\draw (323,41) node [scale=0.7,color={rgb, 255:red, 128; green, 128; blue, 128 }  ,opacity=0.76 ]  {$1< \alpha \ \leq 9$};

	\end{tikzpicture}
	\caption{Set of equilibria in the phase space as $\alpha$ varies.} 
	\label{meq}

\end{figure}

Let us now change slightly the vector $c$ by setting: $c=[\frac{\alpha}{3}, -1, \frac{2\alpha}{3}]'$ with $\alpha \in [0,9]$. Notice that we have multiple equilibria when $\alpha=1 \implies c^*=[\frac 13, -1, \frac 23]'$ as in that case one can check that condition of Theorem \ref{TH:UNIQ}(iv) holds. In Fig. \ref{meq} we show the set of equilibria $\mc X(c)$ in the phase space as $c$ varies as a function of $\alpha$.

Notice that $x^*(c)$ is a piece-wise linear function. We can see that for $0 \le \alpha < 1$ the equilibria (red segment) start from 0, they are unique and located on $\partial \mc L_0^w$, then when $\alpha=1$ (and $ c=c^* $) we have multiple equilibria (orange segment) and finally when $\alpha >1$ the unique equilibria (gray segment) are located on  $\partial \mc L_0^w$ until they eventually reach $w$, which means that all cells output their maximal flow.

We appreciate a phase transition of the dynamical system as the parameter $\alpha$ crosses the value $\alpha=1$. In this case infact, the equilibria undergo a jump discontinuity going from $\underbarr{x}(c^*)$ to $\overbar{x}(c^*)$

\section{Proof of the stability results}\label{sec:analysis}
This section is devoted to prove Theorem 1. We will  first present some technical results concerning properties of the system \eqref{ema} that we will need to prove the main statement. 



We start with the following technical results, whose proofs are presented in Appendix \ref{sec:proof-lemma-invariant+monotone+contraction}, Appendix \ref{sec:proof-minmax-equilibria} and Appendix \ref{sec:proof-lemma-internal-equilibria}  respectively. 

\begin{lemma}\label{LEMMA:INVARIANT+MONOTONE+CONTRACTION}
The dynamical system \eqref{ema} is monotone and non-expansive in $l_1$-distance on $\mathcal{L}_{0}^{w}$.
\end{lemma}

\begin{lemma}\label{LEMMA:MINMAX-EQUILIBRIA}
	The dynamical system \eqref{ema} always admits a maximal equilibrium $\ov x\in\mc L_0^w$ and a minimal equilibrium $\ul x\in\mc L_0^w$. Moreover, the sets \be\label{Xalpha}\mc X_{\alpha}=\l\{x\in\mc L^{\ov x}_{\ul x}:\,\sum_ix_i=\alpha\sum_i\ul x_i+(1-\alpha)\sum_i\ov x_i\r\}\ee for $0\le\alpha\le 1$ are all invariant for \eqref{ema} and, for every initial condition $x(0)\in\mc L_0^w$, the solution of \eqref{ema} is such that $x(t)\stackrel{t\to+\infty}{\longrightarrow}\mc L_{\ul x}^{\ov x}$. 
\end{lemma}

\begin{lemma}\label{LEMMA:INTERNAL-EQUILIBRIA2}
	Let $x^*$ be an equilibrium of the dynamical flow network \eqref{ema} belonging to the interior of the lattice $\mc L_0^w$. Then, there exists an $\varepsilon>0$ such that, every solution of  \eqref{ema} with initial condition $x(0)\in\mc L^w_0$ such that  $||x(0)-x^*||<\eps$, coincides with the solution of the linear dynamics \be\label{linear-system}\dot x=(R'-I)x+c\,. \ee
\end{lemma}


We are now ready to prove a first result that characterizes the set of equilibria.

\begin{proposition}\label{PROP:CURVE-EQUILIBRIA}
	There exists a nondecreasing curve of equilibria joining $\ul{x}$ and $\ov{x}$ with support $\mc X$. 
	Moreover, if $R$ is stochastic irreducible, then such curve is entrywise strictly increasing, while if $R$ is sub-stochastic out-connected such curve is constant so that $\ul{x}=\ov{x}$. 
\end{proposition}

\begin{pf}
For every $0\le\alpha\le1$ the convex compact set $\mc X_{\alpha}$ as defined in \eqref{Xalpha} is invariant for \eqref{ema} by Lemma \ref{LEMMA:MINMAX-EQUILIBRIA}. Then, since $f(x)$ is Lipschtitz-continuous, Lemma 1 in \cite{Lajmanovich.Yorke:1976} implies that \eqref{ema} has at least an equilibrium in $\mc X_{\alpha}$. In fact, observe that, for any $0\le\alpha<\beta\le 1$, if $x^*(\alpha)$ is an equilibrium of \eqref{ema} in $\mc X_{\alpha}$, then the same argument can be applied to show existence of an equilibrium $x^*(\beta)\in\mc X_{\beta}\cap\mc L^{\ov x}_{x^*(\alpha)}$. Moreover, clearly $\lim_{\beta\downarrow\alpha}x^*(\beta)=x^*(\alpha)$. Similarly, one can prove that $\lim_{\beta\uparrow\alpha}x^*(\beta)=x^*(\alpha)$. This shows that there exists a nondecreasing curve of equilibria $[0,1]\ni t\mapsto x^*(t)$ joining $x^*_0=\ul x$ to $x^*_1=\ov x$. 
	
	In order to prove the second part of the claim, fix $0\le\alpha<\beta\le 1$ and let $\mc S\subseteq\{1,\ldots,n\}$ be the set of those cells $i$ such that $x^*_i(\alpha)<x^*_i(\beta)$. 
	If $R$ is stochastic irreducible and $\mc S$ is a strict subset of $\{1,\ldots,n\}$, then
	\be\label{eq:chain}\ba{rcl}
	\beta-\alpha&=&
	\ds\sum_{i\in\mc S}x_i(\beta)-x_i(\alpha)\\
	&=&\ds\sum_{i\in\mc S}S_0^{w_i}\left(\sum\nolimits_jR_{ji}x_j(\beta)+c_i\right)\\
	&&-\ds\sum_{i\in\mc S}S_0^{w_i}\left(\sum\nolimits_jR_{ji}x_j(\alpha)+c_i\right)\\
	&\le&\ds\sum_{i\in\mc S}\sum_{j\in\mc S}R_{ji}(x_j(\beta)-x_j(\alpha))\\
	&<&\ds\sum_{i\in\mc S}(x_j(\beta)-x_j(\alpha))\\
	&=&\beta-\alpha
	\,,\ea\ee
	where the last inequality follows from the fact that $\sum_{i\in\mc S}\sum_{j\in\mc S}R_{ji}z_j<\sum_{i\in\mc S}z_i$ for every positive $z$ and every strict subset $\mc S\subsetneq\{1,\ldots,n\}$. It then follows that necessarily $x_i(\beta)>x_i(\alpha)$ for every $i=1,\ldots,n$. 
	Finally, notice that if $R$ is sub-stochastic out-connected then \eqref{eq:chain} remains valid for every nonempty subset $\mc S\subseteq\{1,\ldots,n\}$, thus implying that necessarily $\ov x=\ul x$ in thus case. 
	\qed
	\end{pf}

We now ready to present the proof of Theorem 1.

\textbf{Proof of Theorem 1}

\begin{enumerate}
	\item[(i)] It immediately follows from Lemma \ref{LEMMA:MINMAX-EQUILIBRIA} and Proposition \ref{PROP:CURVE-EQUILIBRIA} that, when $R$ is sub-stochastic out-connected $\underbarr{x}=\overbar{x}=x^*$ is a global asymptotically stable equilibrium. 
	\item[(ii)]  From Proposition \ref{PROP:CURVE-EQUILIBRIA} we know that there exists a strictly increasing curve joining $\underbarr{x}$ and $\overbar{x}$, which means that either the system has a unique equilibrium or it has a continuum of them. In the latter case, since the curve is strictly increasing,  all non extremal equilibria 
%
 $x^* \in{\mc X}(c)\setminus\{\ul x,\ov x\}$ must belong to the interior of the lattice ${\mc L}_0^w$. All such internal equilibria $x^*$ must satisfy
	\begin{equation}\label{elin3}
	x^*=R'x^*+c 
	\end{equation}
	Now observe that, since $R$ is row-stochastic, we have
	$\1'x=\1' R^{\prime} x+\1' x+\1^{\prime} c \text { so that } \1^{\prime} c=0,$ i.e., for the linear system \eqref{elin3} to admit solutions it is necessary that $c$ is a zero-sum vector. In fact, since the stochastic matrix $R$ is irreducible, we have that $I-R'$ has rank $n-1$  and for every zero-sum vector  $ c $  the set of solutions $x^*$ of the 
	linear system \eqref{elin3} coincides with the line
	\begin{equation}\label{eline}
	\mathcal{H}=\{x^*=Hc+\alpha \pi: \alpha \in \mathbb{R}\}
	\end{equation} 
	Hence, we have that:
	\begin{equation}\label{eq:linequiq}
	\mc X(c)= \mathcal{H} \cap \mathcal{L}_{0}^{w}  = [\underbarr{x}(c), \overbar{x}(c)] 
	\end{equation}
 is a line segment joining $\underbarr{x} \in \partial \mc L_0^w$ and $\overbar{x} \in \partial \mc L_0^w$. 
	
	\item[(iii)] If $\underbarr{x}=\overbar{x}$, then the global convergence follows from Lemma \ref{LEMMA:MINMAX-EQUILIBRIA}. Hence, we need to prove convergence in the case the system admits infinitely many equilibria. 
		Notice that, for every $0\le\alpha\le1$, the set $\mc X_{\alpha}$ defined in \eqref{Xalpha} intersects the line segment $\mc X(c)$ in a single equilibrium point $x^*(\alpha)=\alpha\ul x+(1-\alpha)\ov x$. Moreover, as discussed in the proof of point (ii) above, for every $0<\alpha<1$, such equilibrium $x^*(\alpha)$ belongs to the interior of the lattice ${\mc L_0^w}$, so
	that Lemma \ref{LEMMA:INTERNAL-EQUILIBRIA2} implies that the dynamical flow network \eqref{ema} reduces to the linear dynamical system \eqref{linear-system} in a sufficiently small neighborhood of it. Now observe that all solutions of \eqref{linear-system} with initial condition $x(0)\in\mc X_{\alpha}$ converge to $x^*_{\alpha}$ as $t$ grows large. It then follows that, for every $0\le\alpha\le1$, there exists some $\eps>0$ such that for every solution $x(t)$ of the dynamical flow network with initial condition $x(0)\in\mc X_{\alpha}$ such that $||x(0)-x^*(\alpha)||<\eps$ converges to $x^*(\alpha)$ as $t$ grows large. 	

Now, let $\phi^{t}(x^{\circ})$ be the solution of \eqref{ema} started at $x(0)=x^{\circ}$. 
By Theorem 4.5 in \cite{KhalilBook:02} our last finding implies that, for every $0\le\alpha\le1$ there exists a $\mathcal{K} \mathcal{L} $
	function $\beta(\cdot, \cdot)$  such that $\left\|\phi^{t}(x)-x^{*}(\alpha)\right\| \leq \beta\left(x-x^{*}(\alpha), t\right)$ for every  $x\in\mc X_{\alpha}$ such that $\|x-x^*(\alpha)\|_1\le\eps$. To prove global
	convergence to the set $\mathring{\mc X}(c)$ we need to show that for any $x^{\circ}\in\mc X_{\alpha}$ such that
	$\left\|x^{\circ}-x^{*}\right\|_{1}>\varepsilon,$ there exists a finite time $T \geq 0$ such that $||\phi^{T}\left(x^{\circ}\right)-x^*(\alpha)||_1\le{\varepsilon}$.  For sake of notation, let us put $x^*=x^*(\alpha)$.

	Now	let $ \hat{x}=x^{*}+\frac{\varepsilon}{\left\|x^{\circ}-x^{*}\right\|}\left(x^{\circ}-x^{*}\right),$ for
	which it is easily seen that  $\left\|\hat{x}-x^{*}\right\|_{1}=\varepsilon,$ and
	\begin{align*}\left\|x^{\circ}-x^{*}\right\|_{1}& =\left\|x^{\circ}-\hat{x}\right\|_{1}+\left\|\hat{x}-x^{*}\right\|_{1} \\
	&=\left\|x^{\circ}-\hat{x}\right\|_{1}+\varepsilon \end{align*}
	and consider the trajectories of the system starting from $x^{\circ} $
	and $\hat{x} .$ By the $l_1$-non expansive property ensured by Lemma \ref{LEMMA:INVARIANT+MONOTONE+CONTRACTION} we have
	
	$ \frac{\mathrm{d}}{\mathrm{d} t}\left\|\phi^{t}\left(x^{\circ}\right)-\phi^{t}(\hat{x})\right\|_{1} \leq 0,$ namely		$$\left\|\phi^{t}\left(x^{\circ}\right)-\phi^{t}(\hat{x})\right\|_{1} \leq\left\|x^{\circ}-\hat{x}\right\|_{1}.$$ By the triangle inequality,
	$$
	\begin{aligned}\left\|\phi^{t}\left(x^{\circ}\right)-x^{*}\right\|_{1} & \leq\left\|\phi^{t}\left(x^{\circ}\right)-\phi^{t}(\hat{x})\right\|_{1} \\
	&+\left\|\phi^{t}(\hat{x})-x^{*}\right\|_{1} \\  &=\left\|x^{\circ}-\hat{x}\right\|_{1}+\left\|\phi^{t}(\hat{x})-x^{*}\right\|_{1} \\ &=\left\|x^{\circ}-x^{*}\right\|_{1}-\varepsilon \\
	&+\left\|\phi^{t}(\hat{x})-x^{*}\right\|_{1} \end{aligned}
	$$
	Due to the properties of the $\mathcal{K} \mathcal{L}$ functions, there exists $T_{\frac{\varepsilon}{2}} \geq 0 $
	such that $\beta(x-y, t) \leq \frac{\varepsilon}{2}$ for all $y$ such that $\|y-x^*\|_1\le\eps$ and for all
	$t \geq T_{\frac{\varepsilon}{2}} .$ Thus, we have
	\begin{equation}
	\begin{aligned}\left\|\phi^{t}\left(x^{\circ}\right)-x^{*}\right\|_{1} & \leq\left\|x^{\circ}-x^{*}\right\|_{1}-\varepsilon \\ &
	+\left\|\phi^{t}(\tilde{x})-x^{*}\right\|_{1} \\
	&  \leq\left\|x^{\circ}-x^{*}\right\|_{1}-\frac{\varepsilon}{2} \end{aligned}
	\end{equation}
	for all $t \geq T_{\frac{\varepsilon}{2}} .$ If $||\phi^{T_{\frac{\varepsilon}{2}}}\left(x^{\circ}\right)-x^*(\alpha)||_1\le\eps$ the proof is
	$\text { complete with } T^{-}=T_{\frac{\varepsilon}{2}} $. Otherwise, the same argument
	can be reiterated. Since each step the $\ell_{1}$ distance between
	$\phi^{t}(x) \text { and } x^{*} $ decreases by at least $ \frac{\varepsilon}{2}>0,$ in no more than
	$\left\lceil\frac{2\left\|x^{\circ}-x^{*}\right\|_{1}}{\varepsilon}\right\rceil$ steps, i.e., for $T \leq\left\lceil\frac{2\left\|x^{\circ}-x^{*}\right\|_{1}}{\varepsilon}\right\rceil T_{\frac{\varepsilon}{2}},$ it holds $\|\phi^{T}\left(x^{\circ}\right)-x^*\|_1\le\eps$.  \qed
	\item[(iv)] Because of what said in point (ii) of this proof, the set $\mc X(c)$ has positive length if and only if \eqref{eq:linequiq} defines a non-empty set, i.e. if and only if we can find values of $\alpha \in \R$ such that $0 <Hc+\alpha \pi < w$. Easy computations show that this is the case if and only if 
\begin{equation}\label{key}
	\min _{i}\left\{\frac{(Hc)_{i}}{\pi_{i}}\right\}+\min _{i}\left\{\frac{w_{i}-(Hc)_{i}}{\pi_{i}}\right\}>0  \quad \qed
\end{equation}
\end{enumerate}

\section{Proof of the continuity  results}

This section is devoted to prove Theorem \ref{TH:CONTINUITY}.

We need the following technical results.

\begin{lemma}\label{LEMMA-NONDECMAP} 
	Both $c \mapsto \underline{x}(c)$ and $c \mapsto \bar{x}(c)$ are monotone nondecreasing maps from $\mathbb{R}^{n}$ to $\mathcal{L}_{0}^{w}$.
\end{lemma}
\begin{pf}
	Consider two vectors $c_{1}, c_{2} \in \mathbb{R}^{n}$ such that $c_{1} \leq c_{2}$ and let $x_{1}(t)$ and $x_{2}(t)$
	the two solutions of \eqref{ema} with, respectively, $c=c_{1}$ and $c=c_{2},$
	and with initial condition $x_{1}^\circ=x_{2}^\circ=0 .$ Then we know that these two solutions must converge,
	respectively, to the minimal solutions $\underline{x}\left(c_{1}\right)$ and $\underline{x}\left(c_{2}\right)$. Since the mapping $S_0^w(\cdot)$ is monotone nondecreasing, we can see that
	\begin{align*}
	x_{1}(t) \leq x_{2}(t)  & \Rightarrow \dot{x}_1(t)= S_{0}^{w}\left(R^{\prime} x_{1}(t)+c_{1}\right)-x_1 \\
	& \leq S_{0}^{w}\left(R^{\prime} x_{2}(t)+c_{2}\right)-x_2=\dot{x}_2(t)
	\end{align*}
	since $x_{1}^\circ=x_{2}^\circ=0,$ this implies $x_{1}(t) \leq x_{2}(t)$ for all $t$. This yields $\underline{x}\left(c_{1}\right) \leq \underline{x}\left(c_{2}\right) .$ We
	have proven that $x(c)$ is monotone nondecreasing. The same property for the maximal solution
	$\bar{x}(c)$ follows by an equivalent argument. \qed
	\end{pf}
Lemma \ref{LEMMA-NONDECMAP} allows us to prove the following results that is key for the proof of Theorem \ref{TH:CONTINUITY}:
\begin{lemma}\label{LEMMA-LIM} 
	Let $R\in\R_+^{n\times n}$ be a stochastic matrix and $w\in\R_+^n$ be a nonnegative vector. Then, 
	$$\limsup\limits_{c\to c^*}\bar x(c)=\bar x(c^*)\,,\qquad\liminf\limits_{c\to c^*}\underline x(c)=\underline x(c^*)\,,$$
	for every $c^*\in\R^n$. 
\end{lemma}
\begin{pf}
Let $(c_n)_{n\ge1}$ be any sequence in $\R^n$ such that $c_n\stackrel{n\to+\infty}{\longrightarrow} c^*$ and $\bar x(c_n)\stackrel{n\to+\infty}{\longrightarrow} z^*$. We will show that $z^*\le\ov x(c^*)$. Towards this goal, let $d_n=\sup\{\max\{c_k,c^*\}:\, k\geq n\}$, for $n\ge1$. Clearly, $d_n\stackrel{n\to+\infty}{\longrightarrow}c^*$, while $d_n\geq c_n$, $d_n\ge c^*$, and $d_{n+1}\le d_n$, for every $n\ge1$. Then, Lemma \ref{LEMMA-NONDECMAP} implies that $\bar x(d_n)\geq \bar x(c_n)$, $\bar x(d_n)\ge \bar x(c^*)$, and $\bar x(d_{n+1})\le \bar x(d_{n})$, for every $n\ge1$. Thus, in particular, $\bar x(d_n)$ converges to some $z\in\mc L_0^w$ and such limit satisfies $z\ge z^*$ and $z\ge \ov x(c^*)$. 
On the other hand,  one have that $\ov x(d_n)\in\mc X(d_n)$ is an equilibrium so that $\ov x(d_n)=S_0^w(R'\ov x(d_n)+d_n)$ for every $n\ge1$. By taking the limit of both sides, continuity implies that $z^*=S_0^w(R'z^*+c^*)$ so that $z\in \mc X(c^*)$ must be such that $z\le\ov x(c^*)$. This implies that 
$$\bar x(c^*)=z\ge z^*=\limsup\limits_{c\to c^*}\bar x(c)=\bar x(c^*)\,.$$ 
The liminf part of the statement can then be proven similarly.
\qed	
	\end{pf}

We are now ready to prove Theorem \ref{TH:CONTINUITY}.

\textbf{Proof of Theorem \ref{TH:CONTINUITY}}

\begin{enumerate}
		\item[(i)] Because of Theorem 1, in this case $\mc U=\R^n$ and hence by Lemma \ref{LEMMA-LIM} it follows that, for $c^*\in\R^n$, we have
		\begin{align*}
		\limsup\limits_{c\to c^*}\underline x(c) & \leq \limsup\limits_{c\to c^*}\bar x(c)=\bar x(c^*)=\underline x(c^*) \\
		&=\liminf\limits_{c\to c^*}\underline x(c)\leq \liminf\limits_{c\to c^*}\bar x(c)\
		\end{align*}
		implying that the inequalities in the above must all hold as equalities.
	\item[(ii)]  From  Theorem \ref{TH:UNIQ} it follows that $c^*\in\mc  M$ must satisfy \eqref{eq:nonuniqe}.
	Such condition determines a linear sub-manifold of co-dimension $1$ in $\R^n$. 
		\item[(iii)] It follows from Lemma \ref{LEMMA-LIM} that, for $c^*\in\mc U$, we have
		\begin{align*}
		\limsup\limits_{c\to c^*}\underline x(c) & \leq \limsup\limits_{c\to c^*}\bar x(c)=\bar x(c^*)=\underline x(c^*) \\
		&=\liminf\limits_{c\to c^*}\underline x(c)\leq \liminf\limits_{c\to c^*}\bar x(c)\
		\end{align*}
		implying that the inequalities in the above must all hold as equalities.
			\item[(iv)]  Notice that if $c^*\in\mc M$, then any $c\in\R^n$ such that $c>c^*$ or $c<c^*$ belongs to $\mc U$.  This fact allows one to show that the limit relations in Lemma \ref{LEMMA-LIM} continue to hold true when we restrict $c\in\mc U$ and the proof follows along the same lines. \qed
\end{enumerate}

%


\section{Conclusions and future work}
In this paper we have introduced a nonlinear dynamical system that models a flow dynamic between cells with finite flow capacity. We have completely characterized the set of equilibria of the system and proved the global convergence of the solutions toward this set. Moreover, we have shown how the model exhibits critical phase transitions as the exogenous flow approaches a set of critical values. Future work includes a more in-depth analysis of the discontinuities and their relationship to the network structure and extending the dynamical flow model to allow for nonlinearities in the dependence of the outflow from a cell on the mass of commodity in it.

\bibliography{bibl}

\appendix 
\section{Proof of Lemma \ref{LEMMA:INVARIANT+MONOTONE+CONTRACTION}}\label{sec:proof-lemma-invariant+monotone+contraction}

We first prove that $\mc L_0^w$ is invariant. It is enough to show that when a component $x_i$ reaches the boundary of $\mc L_0^w$, i.e. $x_i=w_i$ or $x_i=0$, then the derivative is non positive or non negative respectively. For $x_i=w_i$, since obviously $  S_0^{w_i}\left( \sum_j R_{ji}x_j+c_i \right) \le w_i$  we have that:
	\begin{equation}\label{key}
\dot{x}_i =S_0^{w_i} \left( \sum_j R_{ji}x_j+c_i \right)-w_i \le 0 
	\end{equation}
	When $x_i=0$,  since $  S_0^{w_i}\left( \sum_j R_{ji}x_j+c_i \right) \ge  \nolinebreak 0$ we have that 
	\begin{equation}\label{key}
\dot{x}_i  = S_0^{w_i} \left( \sum_j R_{ji}x_j+c_i \right)\ge  0 
	\end{equation}
	and this completes the proof. 
	
	We now prove that \eqref{ema} is a monotone system. 
	Set $f_i(x)= S_0^{w_i} \left( \sum_j R_{ji}x_j+c_i \right)-x_i $. It is enough to show that $\displaystyle\frac{\partial f_i}{\partial x_k} \ge 0, \: \forall \: k \neq i $ almost everywhere (i.e. excluding $ 0 $-measure set of points where $f_i$ is not differentiable). 
	
	it is immediate to see that:
	\begin{equation}\label{e4}
	\displaystyle\frac{\partial f_i}{\partial x_k}=
	\begin{cases}
	0 & \mbox{if } \: \sum_j R_{ji}x_j+c_i < 0 \\
	r_{ki} & \mbox{if } \: 0 < \sum_j R_{ji}x_j+c_i < w_i \\
	0 & \mbox{if } \: \sum_j R_{ji}x_j+c_i > w_i \\
	\end{cases} 
	\end{equation}
	Since \eqref{e4} is non negative, therefore (Theorem 1.2) in \cite{Kamke:1929} implies that \eqref{ema} is a monotone system.
	
	Finally, we show that \eqref{ema} is non expansive in $l_1$ distance on $\mc L_0^w$. By monotonicity and using the fact that 
	\begin{equation}
	\sum_{i} \frac{\partial f_{i}}{\partial x_{k}} \leq \sum_{i} r_{k i}-1 \leq 0
	\end{equation}
	the result follows by using (Lemma 5) in \cite{Lovisari.ea:2015}.  \qed

\section{Proof of Lemma \ref{LEMMA:MINMAX-EQUILIBRIA}}\label{sec:proof-minmax-equilibria} From monotonicity and the fact that $\mc L_0^w$ is invariant, the two Cauchy problems
\begin{equation}
\left\{\begin{array}{l}{\dot{x}=S_{0}^{w}\left(R^{\prime} x+c\right)-x} \\ {x_{0}=0}\end{array} \quad\left\{\begin{array}{l}{\dot{x}=S_{0}^{w}\left(R^{\prime} x+c\right)-x} \\ {x_{0}=w}\end{array}\right.\right.
\end{equation}
admit unique solutions that converge to a lower equilibrium $\underline{x}$ and largest equilibrium $\bar{x}$ respectively, i.e. $\underline{x} \leq \bar{x} ;$

Now, let $\ul y=\sum_i\ul x_i$ and $\ov y=\sum_i\ov x_i$. Consider an initial state $x(0)\in\mc L_{\ul x}^{\ov x}$ for $0\le\alpha\le1$. Since the system is non-expansive in $l_1$, both $\|x(t)-\ul x\|_1$ and $\|x(t)-\ov x\|_1$ cannot increase in time, which implies that $\sum_ix_i(t)$ remains constant. It follows that the  sets $\mc X_{\alpha}=\{x\in\mc L^{\ov x}_{\ul x}:\,\sum_ix_i=\alpha \ul y+(1-\alpha)\ov y\}$ are all invariant.

The last claims of the Lemma follow directly from monotonicity. Indeed, for any $x^{\circ}\in\mc L_0^w$, let $\phi^t(x^\circ)$ be the solution of \eqref{ema} at time $t\ge0$. Since $\phi^t(0) \stackrel{t\to+\infty}{\longrightarrow} \underbarr{x}$ and $\phi^t(w) \stackrel{t\to+\infty}{\longrightarrow} \overbar{x}$, then it must be be $\phi^t({x^{\circ}}) \stackrel{t\to+\infty}{\longrightarrow} \mc L_{\underbarr{x}}^{\overbar{x}} \: \: \forall \:  x^{\circ}\in \mc L_0^w$ and in particular, $\forall \: x^{\circ} \in \mc L_{\underbarr{x}}^{\overbar{x}}$. \qed

\section{Proof of Lemma \ref{LEMMA:INTERNAL-EQUILIBRIA2}}\label{sec:proof-lemma-internal-equilibria}
	Observe that an equilibrium $x^* \in$ the interior of ${\mc L_0^w}$  is such that 
$
S_0^w(R'x^*+c)=x^* 
$
belongs to the interior of ${\mc L_0^w}$ 
which in turn implies that 
$$
f(x^*)=S_0^w(R'x^*+c)-x^* = (R'-I)x^*+c\,.  
$$
Since the map $f(x)$ is continuous, there necessarily exists an $\varepsilon>0$ such that for all $||x-x^*||<\eps$ we have that $f(x)= (R'-I)x+c$. 
Since $R$ is stochastic, it has spectral radius in the unitary disk centered in zero so that $R'-I$ has all eigenvalues with nonpositive real part. 
Hence $x^*$ is locally stable (both for the linear dynamical system \eqref{linear-system} and the nonlinear dynamical flow network \eqref{ema}, as they locally coincide), so that we can always find a number $\delta\le \varepsilon $ such that if $\left\|x(0)-x^*\right\|<\delta$ then $\left\|x(t)-x^*\right\|<\varepsilon$ for all $t \ge 0$. This ensures that the trajectories of the system remain in the region where the dynamic is linear and hence the claim follows. \qed

\end{document}